\documentclass[a4paper]{article}
\usepackage{RR}
\usepackage{hyperref}

\usepackage{amssymb}
\usepackage{graphicx}
\usepackage{pstricks}
\usepackage{url}

\RRdate{April 2009}

\RRauthor{Arnaud~Liefooghe \and Laetitia~Jourdan \and El-Ghazali~Talbi}
\authorhead{Liefooghe et al.}
\RRtitle{Un model unifi\'e pour l'optimisation  \'evolutionnaire multiobjectif et son impl\'ementation dans un cadre logiciel g\'en\'erique: ParadisEO-MOEO}
\RRetitle{A Unified Model for Evolutionary Multiobjective Optimization and its Implementation in a General Purpose Software Framework: ParadisEO-MOEO}
\titlehead{A Unified Model for EMO and its Implementation in ParadisEO-MOEO}
\RRresume{
Ce document donne un aperçu concis des algorithmes \'evolutionnaires pour l'optimisation multiobjectif. Un nombre cons\'equent de m\'ethodes \'evolutionnaires d\'edi\'ees \`a la r\'esolution de probl\`emes multiobjectifs a \'et\'e propos\'e \`a ce jour, et une tentative d'unifier les approches existantes est ici pr\'esent\'ee. Sur la base d'une d\'ecomposition \`a grain fin et suite \`a la description des principales questions de {\itshape fitness assignment}, de pr\'eservation de la diversit\'e et d'\'elitisme, un mod\`ele conceptuel est propos\'e et valid\'e en traitant un certain nombre algorithmes classiques comme de simples variantes de la m\^eme structure. Le mod\`ele pr\'esent\'e est alors incorpor\'e dans un logiciel d\'edi\'e \`a la conception et \`a l'impl\'ementation d'algorithmes \'evolutionnaires pour l'optimisation multiobjectif: Paradiseo-MOEO. Cette librairie a d\'emontr\'e sa validit\'e et sa grande souplesse en permettant la r\'esolution d'un grand nombre de probl\`emes d'optimisation multiobjectifs r\'eels et difficiles.
}
\RRabstract{
This paper gives a concise overview of evolutionary algorithms for multiobjective optimization.
A substantial number of evolutionary computation methods for multiobjective problem solving has been proposed so far, 
and an attempt of unifying existing approaches is here presented.
Based on a fine-grained decomposition and following the main issues of fitness assignment, diversity preservation and elitism, 
a conceptual global model is proposed and is validated by regarding a number of state-of-the-art algorithms as simple variants of the same structure.
The presented model is then incorporated into a general-purpose software framework 
dedicated to the design and the implementation of evolutionary multiobjective optimization techniques: ParadisEO-MOEO.
This package has proven its validity and flexibility by enabling the resolution of many real-world and hard multiobjective optimization problems.
}
\RRmotcle{algorithmes \'evolutionnaires, optimisation multiobjectif, mod\`ele conceptuel unifi\'e, conception et impl\'ementation d'algorithme, cadre logiciel}
\RRkeyword{evolutionary algorithms, multiobjective optimization, conceptual unified model, algorithm design and implementation, software frameworks}
\RRprojet{Dolphin}
\RRtheme{\THNum}
\RCLille

\begin{document}
\RRNo{6906}
\makeRR

\section{Introduction}
\label{sec:intro}
Evolutionary Multiobjective Optimization (EMO) is one of the most challenging areas in the field of multicriteria decision making.
Generally speaking, a Multiobjective Optimization Problem (MOP) can be defined by a vector function $f$ of $n\geq2$ objective functions $(f_1,f_2,\dots,f_n)$, 
a set~$X$ of feasible solutions in the \emph{decision space},
and a set~$Z$ of feasible points in the \emph{objective space}.
Without loss of generality, we assume that $Z \subseteq \mathbb{R}^n$ and that all $n$~objective functions are to be minimized. 
To each decision vector $x \in X$ is assigned an objective vector $z \in Z$ on the basis of the vector function~$f : X \rightarrow Z$ with $z = f(x)$. 
A dominance relation is then usually assumed so that a partial order is induced over $X$.
Numerous dominance relations exist in the literature and will be discussed later in the paper.
Let us consider the well-known concept of \emph{Pareto dominance},
for which a given objective vector $z \in Z$ is said to \emph{dominate} another objective vector $z' \in Z$ if 
$\forall i \in \{1,2,\dots,n\}$, $z_i \leq z_i'$ and $\exists j \in \{1,2,\dots,n\} $ such as $z_j < z_j'$. 
An objective vector $z \in Z$ is said to be \emph{nondominated} if there does not exist any other objective vector $z' \in Z$ such that $z'$ dominates $z$.
By extension, we will say that a decision vector $x \in X$ \emph{dominates} a decision vector $x' \in X$ if~$f(x)$ dominates $f(x')$, 
and that a decision vector $x \in X$ is \emph{nondominated} (or \emph{efficient}, \emph{Pareto optimal}) if~$f(x)$ maps to a nondominated point.
The set of all efficient solutions is called \emph{efficient} (or \emph{Pareto optimal}) \emph{set} 
and its mapping in the objective space is called \emph{Pareto front}.

In practice, different resolution scenarios exist and strongly rely on the cooperation between the search process and the decision making process.
Indeed, a distinction can be made between the following forms such a cooperation might take.
For instance, the Decision Maker (DM) may be interested in identifying the whole set of efficient solutions, 
in which case the choice of the most preferred solution is made {\itshape a~posteriori}.
However, when preference information can be provided {\itshape a~priori}, 
the search may lead to the potential best compromise solution(s) over a particular preferred region of the Pareto front.
A third class of methods consists of a progressive, interactive, cooperation between the DM and the solver.
However, in any case, the overall goal is often to identify a set of good-quality solutions.
But generating such a set is usually infeasible, due to the complexity of the underlying problem or to the large number of optima.
Therefore, the overall goal is often to identify a good approximation of it. 
Evolutionary algorithms are commonly used to this end, as they are particularly well-suited to find multiple efficient solutions in a single simulation run.
The reader is referred to~\cite{Deb:01,CLV:07} for more details about EMO.

As pointed out by different authors (see {\itshape e.g.}~\cite{CLV:07,ZLB:04}), approximating an efficient set is itself a bi-objective problem.
Indeed, the approximation to be found must have both good convergence and distribution properties,
as its mapping in the objective space has to be ($i$) close to, and ($ii$) well-spread over the (generally unknown) optimal Pareto front, or a subpart of it.
As a consequence, the main difference between the design of a single-objective and of a multiobjective search method deals with these two goals.
Over the last two decades, major advances, from both algorithmic and theoretical aspect, have been made in the EMO field.
And a large number of algorithms have been proposed.
Among existing approaches, one may cite VEGA~\cite{Sch:85}, MOGA~\cite{FF:93}, NSGA~\cite{SD:94}, NSGA-II~\cite{DA+:02}, 
NPGA~\cite{HNG:94}, SPEA~\cite{ZT:99}, SPEA2~\cite{ZLT:01} or PESA~\cite{CKO:00}.
All these methods are presented and described in~\cite{CLV:07}.
Note that another topic to mention while dealing with EMO relates to performance assessment.
Various quality indicators have been proposed in the literature for evaluating the performance of multiobjective search methods.
The reader is referred to~\cite{ZT+:03} for a review.

In~\cite{ZLB:04}, Zitzler et al. notice that
initial EMO approaches were mainly focused on moving toward the Pareto front~\cite{Sch:85,Fou:85}.
Afterwards, diversity preservation mechanisms quickly emerged~\cite{FF:93,SD:94,HNG:94}.
Then, at the end of the nineteens, the concept of elitism, related to the preservation of nondominated solutions, 
became very popular and is now employed in most recent EMO methods~\cite{ZT:99,ZLT:01,KC:00}.
Specific issues of \emph{fitness assignment}, \emph{diversity preservation} and \emph{elitism} are commonly approved in the community
and are also presented under different names in, for instance,~\cite{CLV:07,ZLB:04}.
Based on these three main notions, several attempts have been made in the past for unifying EMO algorithms.
In~\cite{LZT:00}, the authors focus on elitist EMO search methods.
This study has been later extended in~\cite{ZLB:04} where the algorithmic concepts of fitness assignment, diversity preservation and elitism are largely discussed.
More recently, Deb proposed a robust framework for EMO~\cite{Deb:08} based on NSGA-II (Non-dominated Sorting Genetic Algorithm)~\cite{DA+:02}.
The latter approach is decomposed into three main EMO-components related to elite preservation, nondominated solutions emphasis and diversity maintaining.
However, this model is strictly focused on NSGA-II, whereas other state-of-the-art methods can be decomposed in the same way.
Indeed, a lot of components are shared by many EMO algorithms, so that, 
in somehow, they can all be seen as variants of the same unified model, as it will be highlighted in the remainder of the paper.
Furthermore, some existing models have been used as a basis for the design of tools to help practitioners for MOP solving.
For instance, following~\cite{ZLB:04,LZT:00}, the authors proposed a software framework for EMO called PISA~\cite{BL+:03}.
PISA is a platform and programming language independent interface for search algorithms that consists of two independent modules (the variator and the selector) communicating via text files.
Note that other software frameworks dealing with the design of metaheuristics for EMO have been proposed, 
including jMetal~\cite{DN+:06}, the MOEA toolbox for Matlbax~\cite{TL+:01}, MOMHLib++~\cite{MOMHLib++} and Shark~\cite{Shark}.
These packages will be discussed later in the paper.

The purpose of the present work is twofold.
Firstly, a unified view of EMO is given.
We describe the basic components shared by many algorithms, and we introduce a general purpose model as well as a classification of its fine-grained components.
Next, we confirm its high genericity and modularity by treating a number of state-of-the-art methods as simple instances of the model.
NSGA-II~\cite{DA+:02}, SPEA2~\cite{ZLT:01} and IBEA~\cite{ZK:04} are taken as examples.
Afterwards, we illustrate how this general-purpose model has been used as a starting point for the design and the implementation of 
an open-source software framework dedicated to the reusable design of EMO algorithms, namely ParadisEO-MOEO%
\footnote{ParadisEO-MOEO is available at \url{http://paradiseo.gforge.inria.fr}}.
All the implementation choices have been strongly motivated by the unified view presented in the paper.
This free C++ white-box framework has been widely experimented and has enabled the resolution of a large diversity of MOPs from both academic and real-world applications.
In comparison to the literature, we expect the proposed unified model to be more complete, to provide a more fine-grained decomposition,
and the software framework to offer a more modular implementation than previous similar attempts.
The reminder of the paper is organized as follows.
In Sect.~\ref{sec:model}, a concise, unified and up-to-date presentation of EMO techniques is discussed.
Next, a motivated presentation of the software framework introduced in this paper is given in Sect.~\ref{sec:paradiseo}, 
and is followed by a detailed description of the design and the implementation of EMO algorithms under ParadisEO-MOEO.
Finally, the last section concludes the paper.

\section{The Proposed Unified Model}
\label{sec:model}
An Evolutionary Algorithm (EA)~\cite{ES:03} is a search method that belongs to the class of metaheuristics~\cite{Tal:09}, 
and where a population of solutions is iteratively improved by means of some stochastic operators.
Starting from an initial population, 
each individual is evaluated in the objective space and a selection scheme is performed	to build a so-called parent population.
An offspring population is then created by applying variation operators.
Next, a replacement strategy determines which individuals will survive.
The search process is iterated until a given stopping criterion is satisfied.

As noticed earlier in the paper, in the frame of EMO, 
the main expansions deal with the issues of \emph{fitness assignment}, \emph{diversity preservation} and \emph{elitism}.
Indeed, contrary to single-objective optimization where the fitness value of a solution corresponds to its single objective value in most cases,
a multiobjective fitness assignment scheme is here required to assess the individuals performance, 
as the mapping of a solution in the objective space is now multi-dimensional.
Moreover, trying to approximate the efficient set is not only a question of convergence.
The final approximation also has to be well spread over the objective space, so that a diversity preservation mechanism is usually required.
This fitness and diversity information is necessary to discriminate individuals at the selection and the replacement steps of the EA.
Next, the main purpose of elitism is to avoid the loss of best-found nondominated solutions during the stochastic search process.
These solutions are frequently incorporated into a secondary population, the so-called \emph{archive}.
The update of the archive contents possibly appear at each EA iteration.

As a consequence, whatever the MOP to be solved, the common concepts for the design of an EMO algorithm are the following ones:
\begin {enumerate}
\item Design a representation.
\item Design a population initialization strategy.
\item Design a way of evaluating a solution.
\item Design suitable variation operators. 
\item Decide a fitness assignment strategy.
\item Decide a diversity preservation strategy.
\item Decide a selection strategy.
\item Decide a replacement strategy.
\item Decide an archive management strategy.
\item Decide a continuation strategy.
\end {enumerate}
When dealing with any kind of metaheuristics, one may distinguish problem-specific and generic components.
Indeed, the former four common-concepts presented above strongly depends of the MOP at hand, while the six latter ones can be considered as problem-independent, 
even if some problem-dependent strategies can also be envisaged in some particular cases.
Note that concepts of representation and evaluation are shared by any metaheuristic, 
concepts of population initialization and stopping criterion are shared by any population-based metaheuristic,
concepts of variation operators, selection and replacement are shared by any EA,
whereas concepts of fitness, diversity and archiving are specific to EMO.

\subsection{Components Description}
This section provides a description of components involved in the proposed unified model.
EMO-related components are detailed in more depth.

\subsubsection{Representation}
Solution representation is the starting point for anyone who plans to design any kind of metaheuristic.
A MOP solution needs to be represented both in the decision space and in the objective space.
While the representation in the objective space can be seen as problem-independent,
the representation in the decision space must be relevant to the tackled problem.
Successful applications of metaheuristics strongly requires a proper solution representation.
Various encodings may be used such as binary variables, real-coded vectors, permutations, discrete vectors, and more complex representations.
Note that the choice of a representation will considerably influence the way solutions will be initialized and evaluated in the objective space, 
and the way variation operators will be applied.

\subsubsection{Initialization}
Whatever the algorithmic solution to be designed, a way to initialize a solution (or a population of solutions) is expected.
While dealing with any population-based metaheuristic, 
one has to keep in mind that the initial population must be well diversified in order to prevent a premature convergence.
This remark is even more true for MOPs where the goal is to find a well-converged and a well-spread approximation.
The way to initialize a solution is closely related to the problem under consideration and to the representation at hand.
In most approaches, the initial population is generated randomly or according to a given diversity function.

\subsubsection{Evaluation}
The problem at hand is to optimize a set of objective functions simultaneously over a given search space. 
Then, each time a new solution integrates the population, its objective vector must be evaluated, 
{\itshape i.e.} the value corresponding to each objective function must be set.

\subsubsection{Variation}
The purpose of variation operators is to modify the representation of solutions in order to move them in the search space.
Generally speaking, while dealing with EAs, these problem-dependent operators are stochastic.
Mutation operators are unary operators acting on a single solution whereas recombination (or crossover) operators are mostly binary, and sometimes n-ary. 

\subsubsection{Fitness Assignment}
In the single-objective case, the fitness value assigned to a given solution is most often its unidimensional objective value.
While dealing with MOPs, fitness assignment aims to guide the search toward Pareto optimal solutions for a better convergence.
Extending~\cite{CLV:07,ZLB:04}, we propose to classify existing fitness assignment schemes into four different families:
\begin{itemize}
	\item \emph{Scalar approaches}, where the MOP is reduced to a single-objective optimization problem.
	A popular example consists of combining the $n$ objective functions into a single one by means of a weighted-sum aggregation.
	Other examples are $\epsilon$-constraint or achievement function-based methods~\cite{Mie:99}.
	\item \emph{Criterion-based approaches}, where each objective function is treated separately.
	For instance, in VEGA (Vector Evaluated GA)~\cite{Sch:85}, 
	a parallel selection is performed where solutions are discerned according to their values on a single objective function, independently to the others.
	In lexicographic methods~\cite{Fou:85}, a hierarchical order is defined between objective functions.
	\item \emph{Dominance-based approaches}, where a dominance relation is used to classify solutions.
	For instance, \emph{dominance-rank} techniques computes the number of population items that dominate a given solution~\cite{FF:93}.
	Such a strategy take part in, {\itshape e.g.}, Fonseca and Fleming MOGA (Multiobjective GA)~\cite{FF:93}.
	In \emph{dominance-count} techniques, the fitness value of a solution corresponds to the number of individuals that are dominated by these solutions~\cite{ZT:99}.
	Finally, \emph{dominance-depth} strategies consists of classifying a set of solutions into different classes (or fronts)~\cite{Gol:89}.
	Hence, a solution that belongs to a class does not dominate another one from the same class;
	so that individuals from the first front all belong to the best nondominated set, 
	individuals from the second front all belong to the second best nondominated set, and so on.
	The latter approach is used in NSGA (Non-dominated Sorting GA)~\cite{SD:94} and NSGA-II~\cite{DA+:02}.
	However, note that several schemes can also be combined, what is the case, for example, in~\cite{ZT:99}.
	In the frame of dominance-based approaches, the most commonly used dominance relation is based on Pareto-dominance as given in Sect.~\ref{sec:intro}.
	But some recent techniques are based on other dominance operators such as $\epsilon$-dominance in~\cite{DMM:05} or g-dominance in~\cite{MS+:08}.
	\item \emph{Indicator-based approaches}, where the fitness values are computed by comparing individuals on the basis of a quality indicator $I$.
	The chosen indicator represents the overall goal of the search process.
	Generally speaking, no particular diversity preservation mechanism is in usual necessary, with regards to the indicator being used.
	Examples of indicator-based EAs are IBEA (Indicator-Based EA)~\cite{ZK:04} or SMS-EMOA (S-Metric Selection EMO Algorithm)~\cite{BNE:07}.
\end{itemize}

\subsubsection{Diversity Assignment}
As noticed in the previous section, aiming at approximating the efficient set is not only a question of convergence.
The final approximation also has to be well spread over the objective space.
However, classical dominance-based fitness assignment schemes often tend to produce premature convergence by privileging nondominated solutions, 
what does not guarantee a uniformly sampled output set.
In order to prevent that issue, a diversity preservation mechanism, based on a given distance measure, is usually integrated into the algorithm 
to uniformly distribute the population over the trade-off surface.
In the frame of EMO, a common distance measure is based on the euclidean distance between objective vectors.
But, this measure can also be defined in the decision space or can even combined both spaces.
Popular examples of EMO diversity assignment techniques are sharing or crowding.
%
The notion of \emph{sharing} (or \emph{fitness sharing}) has initially been suggested by Goldberg and Richardson~\cite{GR:87} to preserve diversity among the solutions of an EA population.
It has first been employed by Fonseca and Fleming~\cite{FF:93} in the frame of EMO.
This \emph{kernel} method consists of estimating the distribution density of a solution using a so-called \emph{sharing function} 
that is related to the sum of distances to its neighborhood solutions.
A sharing distance parameter specifies the similarity threshold, {\itshape i.e.} the size of \emph{niches}.
The distance measure between two solutions can be defined in the decision space, in the objective space or can even combined both. 
Nevertheless, a distance metric partly or fully defined in the parameter space strongly depends of the tackled problem. 
Another diversity assignment scheme is the concept of \emph{crowding}, 
firstly suggested by Holland~\cite{Hol:75} and used by De~Jong to prevent \emph{genetic drift}~\cite{Dej:75}.
It is employed by Deb et al.~\cite{DA+:02} in the frame of NSGA-II.
Contrary to sharing, this scheme allows to maintain diversity without specifying any parameter.
It consists in estimating the density of solutions surrounding a particular point of the objective space.

\subsubsection{Selection}
\label{sec:select}
The selection step is one of the main search operators of EAs.
It consists of choosing some solutions that will be used to generate the offspring population.
In general, the better is an individual, the higher is its chance of being selected.
Common strategies are deterministic or stochastic tournament, roulette-wheel selection, random selection, etc. 
An existing EMO-specific elitist scheme consists of including solutions from the archive in the selection process, 
so that nondominated solutions also contributes to the evolution engine.
Such an approach has successfully been applied in various elitist EMO algorithms including SPEA~\cite{ZT:99}, SPEA2~\cite{ZLT:01} or PESA~\cite{CKO:00}.
In addition, in order to prohibit the crossover of dissimilar parents, 
mating restriction~\cite{Gol:89} can also be mentioned as a candidate strategy to be integrated into EMO algorithms.

\subsubsection{Replacement}
Selection pressure is also affected at the replacement step where survivors are selected from both the current and the offspring population.
In generational replacement, the offspring population systematically replace the parent one.
An elitist strategy consists of selecting the $N$ best solutions from both populations, where $N$ stands for the appropriate population size.

\subsubsection{Elitism}
Another essential issue about MOP solving is the notion of \emph{elitism}.
It mainly consists of maintaining an external set, the so-called \emph{archive}, 
that allows to store either all or a subset of nondominated solutions found during the search process.
This secondary population mainly aims at preventing the loss of these solutions during the stochastic optimization process. 
The update of the archive contents with new potential nondominated solutions is mostly based on the Pareto-dominance criteria.
But, in the literature, other dominance criterion are found and can be used instead of the Pareto-dominance relation.
Examples are weak-dominance, strict-dominance, $\epsilon$-dominance~\cite{HP:94}, etc.
When dealing about archiving, one may distinguished four different techniques depending on the problem properties, the designed algorithm and the number of desired solutions:
($i$) \emph{no archive}, ($ii$) an \emph{unbounded archive}, ($iii$)~a \emph{bounded archive} or ($iv$) a \emph{fixed-size archive}.
Firstly, if the current approximation is maintained by, or contained in the main population itself, there can be no archive at all.
On the other hand, if an archive is maintained, it usually comprises the current nondominated set approximation, as dominated solutions are removed.
Then, an unbounded archive can be used in order to save the whole set of nondominated solutions found until the beginning of the search process.
However, as some continuous optimization problems may contain an infinite number of nondominated solutions, it is simply not possible to save them all.
Therefore, additional operations must be used to reduce the number of stored solutions.
Then, a common strategy is to bound the size of the archive according to some fitness and/or diversity assignment scheme(s).
Finally, another archiving technique consists of a fixed size storage capacity, where a bounding mechanism is used when there is too much nondominated solutions, 
and some dominated solutions are integrated in the archive if the nondominated set is too small, what is done for instance in SPEA2~\cite{ZLT:01}.
Usually, an archive is used as an external storage only.
However, archive members can also be integrated during the selection phase of an EMO algorithm~\cite{ZT:99}, see Sect.~\ref{sec:select}.

\subsubsection{Stopping criteria}
Since an iterative method computes successive approximations, a practical test is required to determine when the process must stop.
Popular examples are  a given number of iterations, a given number of evaluations, a given run time, etc.

\subsection{State-of-the-art EMO Methods as Instances of the Proposed Model}
\label{sec:instances}
By means of the unified model proposed in this paper, we claim that a large number of state-of-the-art EMO algorithms proposed in the last two decades 
are based on variations of the problem-independent components presented above.
In Table~\ref{tab:instances}, three EMO approaches, namely NSGA-II~\cite{DA+:02}, SPEA2~\cite{ZLT:01} and IBEA~\cite{ZK:04}, 
are regarded as simple instances of the unified model proposed in this paper.
Of course, only problem-independent components are presented.
NSGA-II and SPEA2 are two of the most frequently encountered EMO algorithms of the literature, 
either for tackling an original MOP or to serve as references for comparison.
Regarding IBEA, it is a good illustration of the new EMO trend dealing with indicator-based search that started to become popular in recent years.
We can see in the table that these three state-of-the-art algorithms perfectly fit into our unified model for EMO, what strongly validates the proposed approach.
But other examples can be found in the literature.
For instance, the only components that differ from NSGA~\cite{SD:94} to NSGA-II~\cite{DA+:02} is the diversity preservation strategy, 
that is based on sharing in NSGA and on crowding in NSGA-II.
Another example is the $\epsilon$-MOEA proposed in~\cite{DMM:05}.
This algorithm is a modified version of NSGA-II where the Pareto-dominance relation used for fitness assignment has been replaced by the $\epsilon$-dominance relation.
Similarly, the g-dominance relation proposed in~\cite{MS+:08} is experimented by the authors on a NSGA-II-like EMO technique 
where the dominance relation has been modified in order to take the DM preferences into account by means of a reference point.
\begin{table*}
\caption{State-of-the-art EMO methods as instances of the proposed unified model.}
\begin{center}
\rotatebox{90}{
\begin{tabular}{|l|c|c|c|}
\hline
\multicolumn{1}{|c|}{\emph{Components}}  
&  \multicolumn{1}{|c|}{\emph{NSGA-II}~\cite{DA+:02}}
&  \multicolumn{1}{|c|}{\emph{SPEA2}~\cite{ZLT:01}}
&  \multicolumn{1}{|c|}{\emph{IBEA}~\cite{ZK:04}}  \\
\hline
\emph{Fitness assignment}    &  dominance-depth      &  dominance-count and rank  &  binary quality indicator  \\
\hline
\emph{Diversity assignment}  &  crowding distance    &  nearest neighbor          & none  \\
\hline
\emph{Selection}             &  binary tournament    &  binary  tournament        &  binary  tournament  \\
\hline
\emph{Replacement}           &  elitist replacement  &  generational replacement  &  elitist replacement  \\
\hline
\emph{Archiving}             &  none                 &  fixed-size archive        &  none  \\
\hline
\emph{Stopping criteria}     & max. number of generations  &  max. number of generations  &  max. number of generations  \\
\hline
\end{tabular}
}
\label{tab:instances}
\end{center}
\end{table*}
%

\section{Design and Implementation under ParadisEO-MOEO}
\label{sec:paradiseo}
In this section, we provide a general presentation of ParadisEO, a software framework dedicated to the design of metaheuristics,
and a detailed description of the ParadisEO module specifically dedicated to EMO, namely ParadisEO-MOEO.
Historically, ParadisEO was especially dedicated to parallel and distributed metaheuristics 
and was the result of the PhD work of S\'ebastien Cahon, supervised by Nouredine Melab and El-Ghazali Talbi~\cite{CMT:04}.
The initial version already contained a few number of EMO-related features, mainly with regard to archiving.
This work has been partially extended and presented in~\cite{LB+:07}.
But since then, the ParadisEO-MOEO module has been completely redesigned in order to confer an even more fine-grained decomposition 
in accordance with the unified model presented above.

\subsection{Motivations}
In practice, there exists a large diversity of optimization problems to be solved, 
engendering a wide number of possible models to handle in the frame of a metaheuristic solution method.
Moreover, a growing number of general-purpose search methods are proposed in the literature, with evolving complex mechanisms.
From a practitioner point of view, there is a popular demand to provide a set of ready-to-use metaheuristic implementations, allowing a minimum programming effort.
On the other hand, an expert generally wants to be able to design new algorithms, to integrate new components into an existing method, 
or even to combine different search mechanisms.
As a consequence, an approved approach for the development of metaheuristics is the use of frameworks.
A metaheuristic software framework may be defined by a set of components based on a strong conceptual separation of the invariant part and the problem-specific part of metaheuristics.
Then, each time a new optimization problem is tackled, both code and design can directly be reused in order to redo as little code as possible.

\subsection{ParadisEO and ParadisEO-MOEO}
ParadisEO\footnote{\url{http://paradiseo.gforge.inria.fr}} is a white-box object-oriented software framework dedicated to the flexible design of metaheuristics 
for optimization problems of both discrete and combinatorial nature.
Based on EO (Evolving Objects)\footnote{\url{http://eodev.sourceforge.net}}~\cite{KM+:01}, 
this template-based, ANSI-C++ compliant computation library is portable across both Unix-like and Windows systems.
Moreover, it tends to be used both by non-specialists and optimization experts.
ParadisEO is composed of four connected modules that constitute a global framework.
Each module is based on a clear conceptual separation of the solution methods from the problems they are intended to solve.
This separation confers a maximum code and design reuse to the user.
The first module, ParadisEO-EO, provides a broad range of components for the development of population-based metaheuristics, including evolutionary algorithms or particle swarm optimization techniques.
Second, ParadisEO-MO contains a set of tools for single-solution based metaheuristics, {\itshape i.e.} local search, simulated annealing, tabu search, etc.
Next, ParadisEO-MOEO is specifically dedicated to the reusable design of metaheuristics for multiobjective optimization.
Finally, ParadisEO-PEO provides a powerful set of classes for the design of parallel and distributed metaheuristics: 
parallel evaluation of solutions, parallel evaluation function, island model and cellular model.
In the frame of this paper, we will exclusively focus on the module devoted to multiobjective optimization, namely ParadisEO-MOEO.

ParadisEO-MOEO provides a flexible and modular framework for the design of EMO metaheuristics.
Its implementation is based on the unified model proposed in the previous section and is conceptually divided into fine-grained components.
On each level of its architecture, a set of abstract classes is proposed and a wide range of instantiable classes, 
corresponding to different state-of-the-art strategies,  are also provided.
Moreover, as the framework aims to be extensible, flexible and easily adaptable, 
all its components are generic so that its modular architecture allows to quickly and conveniently develop any new scheme with a minimum code writing. 
The underlying goal here is to follow new strategies coming from the literature and, if need be, to provide any additional components required for their implementation.
ParadisEO-MOEO constantly evolves and new features might be regularly added to the framework in order to provide a wide range of efficient and modern concepts  
and to reflect the most recent advances of the EMO field.

\subsection{Main Characteristics}
A framework is usually intended to be exploited by a large number of users.
Its exploitation could only be successful if a range of user criteria are satisfied.
Therefore, the main goals of the ParadisEO software framework are the following ones:
\begin{itemize}
\item{\emph{Maximum design and code reuse.}}
The framework must provide a whole architecture design for the metaheuristic approach to be used.
Moreover, the programmer may redo as little code as possible.
This aim requires a clear and maximal conceptual separation of the solution methods and the problem to be solved.
The user might only write the minimal problem-specific code and the development process might be done in an incremental way,
what will considerably simplify the implementation and reduce the development time and cost.
\item{\emph{Flexibility and adaptability.}}
It must be possible to easily add new features or to modify existing ones without involving other components.
Users must have access to source code and use inheritance or specialization concepts of object-oriented programming to derive new components from base or abstract classes.
Furthermore, as existing problems evolve and new others arise, the framework components must be conveniently specialized and adapted.
\item{\emph{Utility.}}
The framework must cover a broad range of metaheuristics, fine-grained components, problems, parallel and distributed models, hybridization mechanisms, etc.
\item{\emph{Transparent and easy access to performance and robustness.}}
As the optimization applications are often time-consuming, the performance issue is crucial.
Parallelism and distribution are two important ways to achieve high performance execution.
Moreover, the execution of the algorithms must be robust in order to guarantee the reliability and the quality of the results.
Hybridization mechanisms generally allow to obtain robust and better solutions.
\item{\emph{Portability.}}
In order to satisfy a large number of users, the framework must support many material architectures (sequential, parallel, distributed) 
and their associated operating systems (Windows, Linux, MacOS).
\item{\emph{Usability and efficiency.}}
The framework must be easy to use and must not contain any additional cost in terms of time or space complexity 
in order to keep the efficiency of a special-purpose implementation.
On the contrary, the framework is intented to be less error-prone than a specifically developed metaheuristic.
\end{itemize}

The ParadisEO platform honors all the above-mentioned criteria and aims to be used by both non-specialists and optimization experts.
Furthermore, The ParadisEO-MOEO module must cover additional goals related to EMO.
Thus, in terms of design, it might for instance be a commonplace to extend a single-objective optimization problem to the multiobjective case 
without modifying the whole metaheuristic implementation.

\subsection{Existing Software Frameworks for Evolutionary Multiobjective Optimization}
Many frameworks dedicated to the design of metaheuristics have been proposed so far. 
However, very few are able to handle MOPs, even if some of them provide components for a few particular EMO strategies, 
such as ECJ~\cite{ECJ}, JavaEVA~\cite{SU:05} or Open~BEAGLE~\cite{GP:06}.
Table~\ref{tab:frameworks} gives a non-exhaustive comparison between a number of existing software frameworks for EMO, 
including jMetal~\cite{DN+:06}, the MOEA toolbox for Matlab a~\cite{TL+:01}, MOMHLib++~\cite{MOMHLib++}, PISA~\cite{BL+:03} and Shark~\cite{Shark}.
Note that other software packages exist for multiobjective optimization~\cite{PVS:08}, but some cannot be considered as frameworks and others do not deal with EMO.
The frameworks presented in Table~\ref{tab:frameworks} are distinguished according to the following criteria: 
the kind of MOPs they are able to tackle (continuous and/or combinatorial problems), 
the availability of statistical tools (including performance metrics), 
the availability of hybridization or parallel features, 
the framework type (black box or white box), the programming language and the license type (free or commercial).
\begin{table*}
\caption{Main characteristics of some existing frameworks for multiobjective metaheuristics.}
\label{tab:frameworks}
\centering
\rotatebox{90}{
\begin{tabular}{|c|c|c|c|c|c|c|c|c|c|}
\hline
Framework       &  \multicolumn{2}{c|}{Problems}  &  \multicolumn{2}{c|}{Statistical tools}  &  Hybrid.  &  Parallel  &  Type   &  Lang.  &  License     \\
                &  Cont.  &  Comb.                &  Off-line  &  On-line                    &           &            &         &         &              \\
\hline
jMetal          &  yes         &  yes             &  yes       &  no                         &  yes      &  no        &  white  &  java   &  free        \\
\hline
MOEA for Matlab &  yes         &  no              &  no        &  no                         &  no       &  yes       &  black  &  matlab &  free / com. \\
\hline
MOMHLib++       &  yes         &  yes             &  no        &  no                         &  yes      &  no        &  white  &  c++    &  free        \\
\hline
PISA            &  yes         &  yes             &  yes       &  no                         &  no       &  no        &  black  &  any    &  free        \\
\hline
Shark           &  yes         &  no              &  no        &  no                         &  yes      &  no        &  white  &  c++    &  free        \\
\hline
\hline
ParadisEO       &  yes         &  yes             &  yes       &  yes                        &  yes      &  yes       &  white  &  c++    &  free        \\
\hline
\end{tabular}
}
\end{table*}
Firstly, let us mentioned that every listed software framework is free of use, except the MOEA toolbox designed for the commercial-software Matlab.
They can all handle continuous problem, but only a subpart is able to deal with combinatorial MOPs.
Moreover, some cannot be considered as white-box frameworks since their architecture is not decomposed into components.
For instance, to design a new algorithm under PISA, it is necessary to implement it from scratch, as no existing element can be reused.
Similarly, even if Shark can be considered as a white-box framework, its components are not as fine-grained as the ones of ParadisEO.
On the contrary, ParadisEO is an open platform where anyone can contribute and add his/her own features.
Finally, only a few ones are able to deal with hybrid and parallel metaheuristics at the same time.
Hence, with regards to the taxonomy proposed in~\cite{Tal:02}, only relay hybrid metaheuristics can be easily implemented within jMetal, MOMHLib++ and Shark,
whereas ParadisEO provides tools for the design of all classes of hybrid models, including teamwork hybridization.
Furthermore, in opposition to jMetal and MOMHLib++, ParadisEO offers easy-to-use models for the design of parallel and distributed EMO algorithms.
Therefore, ParadisEO seems to be the only existing software framework that achieves all the aforementioned goals.

\subsection{Implementation}
This section gives a detailed description of the base classes provided within the ParadisEO framework to design an EMO algorithm%
\footnote{The classes presented in this paper are described as in version $1.2$ of ParadisEO.}.
The flexibility of the framework and its modular architecture based on the three main multiobjective metaheuristic design issues 
(fitness assignment, diversity preservation and elitism) allows to implement efficient algorithms in solving a large diversity of MOPs.
The granular decomposition of ParadisEO-MOEO is based on the unified model proposed in the previous section.

As an EMO algorithm differs of a mono-objective one only in a number of points, some ParadisEO-EO components are directly reusable in the frame of ParadisEO-MOEO.
Therefore, in the following, note that the names of ParadisEO-EO classes are all prefixed by \verb|eo| 
whereas the names of ParadisEO-MOEO classes are prefixed by \verb|moeo|.
ParadisEO is an object-oriented platform, so that its components will be specified by the UML standard~\cite{UML}.
But, due to space limitation, only a subpart of UML diagrams are given, 
but the whole inheritance diagram as well as class documentation and many examples of use are available on the ParadisEO website.
Moreover, a large part of ParadisEO components are based on the notion of \emph{template} and are defined as class templates.
This concept and many related functions are featured within the C++ programming language and allows classes to handle generic types,
so that they can work with many different data types without having to be rewritten for each one.

In the following, both problem-dependent and problem-independent components are detailed.
Hence, basic (representation, evaluation, initialization and stopping criteria), EMO-specific (fitness, diversity and elitism) and EA-related (variation, selection, replacement) components are outlined.
Finally, the way to build a whole EMO algorithm is presented and a brief discussion concludes the section.

\subsubsection{Representation}
\label{sec:representation}
\begin{figure}[htbp]
	\centering
		\includegraphics[scale=0.24]{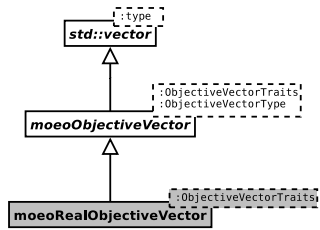}
	\caption{UML diagram for the representation of a solution in the objective space.}
	\label{fig:uml:objectivevector}
\end{figure}

\begin{figure}[htbp]
	\centering
		\includegraphics[scale=0.24]{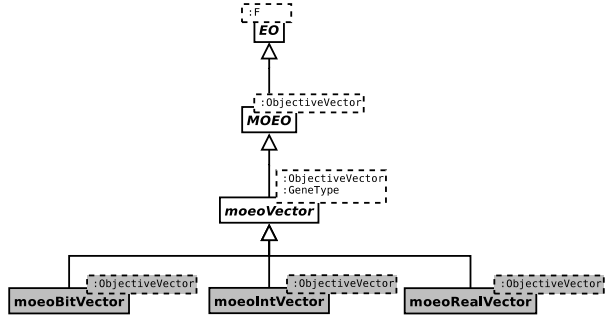}
	\caption{UML diagram for the representation of a solution.}
	\label{fig:uml:moeo}
\end{figure}

A solution needs to be represented both in the decision space and in the objective space.
While the representation in the objective space can be seen as problem-independent, the representation in the decision space must be relevant to the tackled problem.
Using ParadisEO-MOEO, the first thing to do is to set the number of objectives for the problem under consideration and, for each one, 
if it has to be minimized or maximized.
Then, a class inheriting of \verb|moeoObjectiveVector| has to be created for the representation of an objective vector, as illustrated in Fig.~\ref{fig:uml:objectivevector}.
Besides, as a big majority of MOPs deals with real-coded objective values, a class modeling real-coded objective vectors is already provided.
Note that this class can also be used for any MOP without loss of generality.
Next, the class used to represent a solution within ParadisEO-MOEO must extend the \verb|MOEO| class in order to be used for a specific problem.
This modeling tends to be applicable for every kind of problem with the aim of being as general as possible.
Nevertheless, ParadisEO-MOEO also provides easy-to-use classes for standard vector-based representations and, in particular, 
implementations for vectors composed of bits, of integers or of real-coded values that can thus directly be used in a ParadisEO-MOEO-designed application. 
These classes are summarized in Fig.~\ref{fig:uml:moeo}.

\subsubsection{Initialization}
A number of initialization schemes already exists in a lot of libraries for standard representations, what is also the case within ParadisEO.
But some situations could require a combination of many operators or a specific implementation.
Indeed, the framework provides a range of initializers all inheriting of \verb|eoInit|, as well as an easy way to combine them thanks to a \verb|eoCombinedInit| object.

\subsubsection{Evaluation}
The way to evaluate a given solution must be ensured by components inheriting of the \verb|eoEvalFunc| abstract class.
It basically takes a \verb|MOEO| object and sets its objective vector.
Generally speaking, for real-world optimization problems, evaluating a solution in the objective space is by far the most computationally expensive step of any metaheuristic.
A possible way to overcome this trouble is the use of parallel and distributed models,
that can largely be simplify in the frame of ParadisEO thanks to the ParadisEO-PEO module of the software library.
The reader is referred to~\cite{CMT:04} for more information on how to parallelize the evaluation step of a metaheuristic within ParadisEO-PEO.

\subsubsection{Variation}
All variation operators must derive from the \verb|eoOp| base class.
Four abstract classes inherit of \verb|eoOp|, namely \verb|eoMonOp| for mutation operators,
\verb|eoBinOp| and \verb|eoQuadOp| for recombination operators and \verb|eoGenOp| for other kinds of variation operators.
Various operators of same arity can also be combined using some helper classes.
Note that variation mechanisms for some classical (real-coded, vector-based or permutation-based) representations are already provided in the framework.
Moreover, an hybrid mechanism can easily be designed by using a mono-objective local search as mutation operator, as they both inherit of the same class,
see ParadisEO-MO~\cite{BJT:08}.
The set of all variation operators designed for a given problem must be embedded into a \verb|eoTranform| object.

\subsubsection{Fitness Assignment}
\label{sec:fitness}
\begin{figure}[htbp]
	\centering
		\includegraphics[scale=0.24]{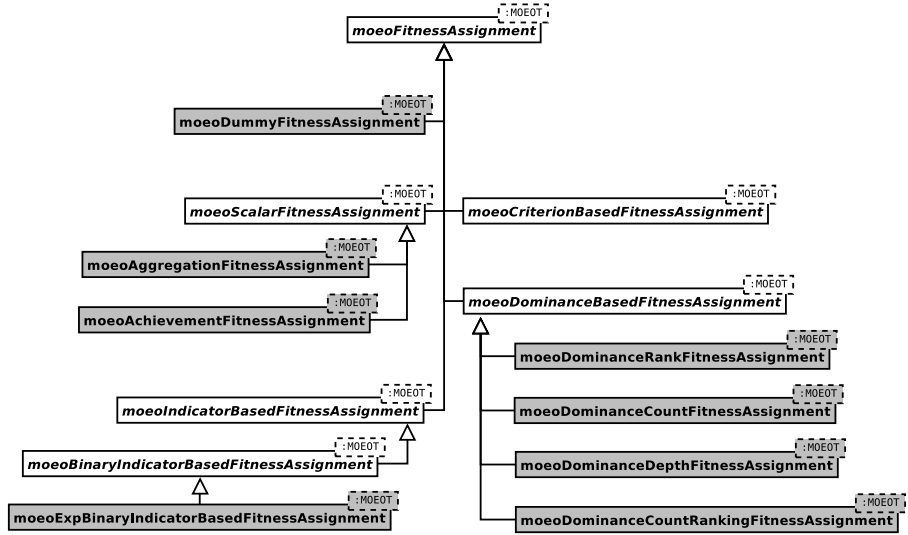}
	\caption{UML diagram for fitness assignment.}
	\label{fig:uml:fitness}
\end{figure}

Following the taxonomy introduced in Sect.~\ref{sec:model}, the fitness assignment schemes are classified into four main categories, 
as illustrated in the UML diagram of Fig.~\ref{fig:uml:fitness}:
scalar approaches, criterion-based approaches, dominance-based approaches and indicator-based approaches.
Non-abstract fitness assignment schemes provided within ParadisEO-MOEO are 
the \emph{achievement scalarizing functions}, the \emph{dominance-rank}, \emph{dominance-count} and \emph{dominance-depth} schemes, 
as well as the \emph{indicator-based} fitness assignment strategy proposed in~\cite{ZK:04}.
Moreover, a dummy fitness assignment strategy has been added in case it would be useful for some specific implementation.

\subsubsection{Diversity Assignment}
\begin{figure}[htbp]
	\centering
		\includegraphics[scale=0.24]{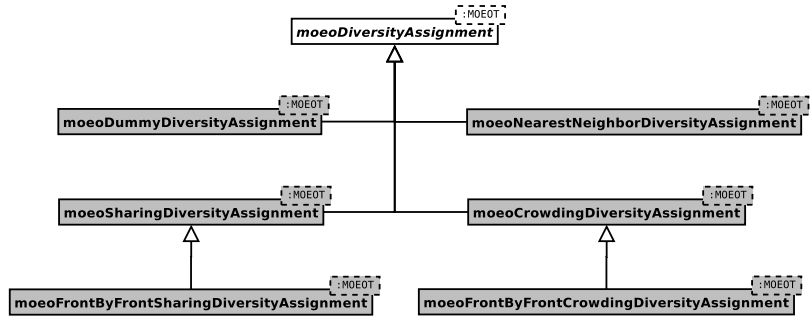}
	\caption{UML diagram for diversity assignment.}
	\label{fig:uml:diversity}
\end{figure}

As illustrated in Fig.~\ref{fig:uml:diversity}, the diversity preservation strategy to be used must inherit of the \verb|moeoDiversityAssignment| class.
Hence, in addition to a dummy technique, a number diversity assignment schemes are already available, 
including sharing~\cite{GR:87}, crowding~\cite{Hol:75} and a nearest neighbor scheme~\cite{ZLT:01}.

\subsubsection{Selection}
There exists a large number of selection strategies in the frame of EMO.
Four ones are provided within ParadisEO-MOEO.
First, the \emph{random selection} consists of selecting a parent randomly among the population members, without taking nor the fitness or the diversity information into account.
Second, the \emph{deterministic tournament selection} consists of performing a tournament between $m$ randomly chosen population members and in selecting the best one.
Next, the \emph{stochastic tournament selection} consists of performing a binary tournament between randomly chosen population members 
and in selecting the best one with a probability $p$ or the worst one with a probability~$(1-p)$.
Finally, the \emph{elitist selection} consists of selecting a population member based on some selection scheme with a probability $p$, 
or in selecting an archive member using another selection scheme with a probability $(1-p)$.
Thus, nondominated (or most-preferred) solutions also contribute to the evolution engine by being used as parents.
A selection method 
needs to be embedded into a \verb|eoSelect| object to be properly used.
Of course, everything is done to easily implement a new selection scheme with a minimum programming effort.

\subsubsection{Replacement}
A large majority of replacement strategies depends on the fitness and/or the diversity value(s) and can then be seen as EMO-specific.
Three replacement schemes are provided within ParadisEO-MOEO, but this list is not exhaustive as new ones can easily be implemented due to the genericity of the framework.
First, the \emph{generational replacement} consists of keeping the offspring population only, while all parents are deleted.
Next, the \emph{one-shot elitist replacement} consists of preserving the $N$ best solutions, where $N$ stands for the population size.
At last, the \emph{iterative elitist replacement} consists of repeatedly removing the worst solution until the required population size is reached. 
Fitness and diversity information of remaining individuals is updated each time there is a deletion.

\subsubsection{Elitism}
\begin{figure}[htbp]
	\centering
		\includegraphics[scale=0.24]{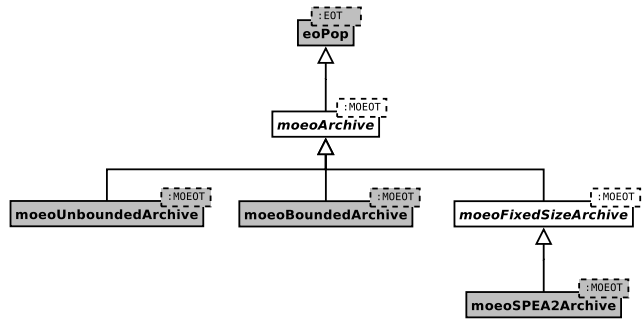}
	\caption{UML diagram for archiving.}
	\label{fig:uml:archive}
\end{figure}

\begin{figure}[htbp]
	\centering
		\includegraphics[scale=0.24]{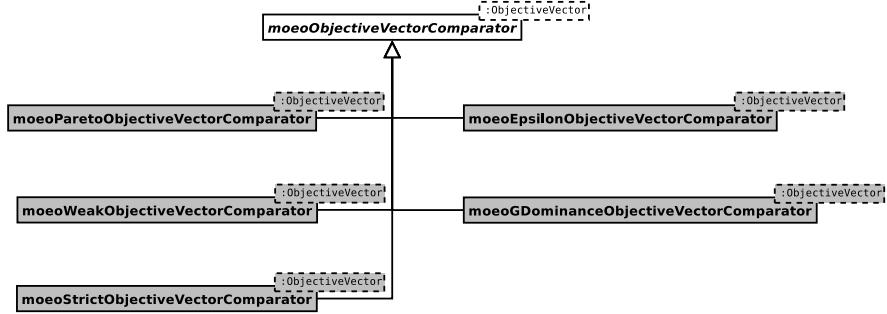}
	\caption{UML diagram for for dominance relation (used for pairwise objective vector comparison).}
	\label{fig:uml:comparator.objectivevector}
\end{figure}

As shown in Fig.~\ref{fig:uml:archive}, in terms of implementation, an archive is represented by the \verb|moeoArchive| abstract class
and is a population using a particular dominance relation to update its contents.
An abstract class for fixed-size archives is given, but implementations of an unbounded archive, 
a general-purpose bounded archive based on a fitness and/or a diversity assignment scheme(s) as well as the SPEA2 archive are provided.
Furthermore, as shown in Fig.~\ref{fig:uml:comparator.objectivevector}, 
ParadisEO-MOEO offers the opportunity to use different dominance relation to update an archive contents by means of a \verb|moeoObjectiveVectorComparator| object, 
including Pareto-dominance, weak-dominance, strict-dominance, $\epsilon$-dominance~\cite{HP:94}, and g-dominance~\cite{MS+:08}.

\subsubsection{Stopping Criteria, Checkpointing and Statistical Tools}
\label{sec:stat}
In the frame of ParadisEO, many stopping criteria extending \verb|eoContinue| are provided.
For instance, the algorithm can stop after a given number of iterations, a given number of evaluations,
a given run time or in an interactive way as soon as the user decides to.
Moreover, different stopping criterion can be combined, in which case the process stops once one of the embedded criteria is satisfied.
In addition, many other procedures may be called at each iteration of the main algorithm.
The \verb|eoCheckPoint| class allows to perform some systematic actions at each algorithm iteration in a transparent way 
by being integrated into the global \verb|eoContinue| object.
The checkpointing engine is particularly helpful for fault tolerance mechanisms and to compute statistical tools.
Indeed, some statistical tools are also provided within ParadisEO-MOEO.
Then, it is for instance possible to save the contents of the current approximation set at each iteration, 
so that the evolution of the current nondominated front can be observed or study using graphical tools such as GUIMOO%
\footnote{GUIMOO is a Graphical User Interface for Multiobjective Optimization available at {\tt http://guimoo.gforge.inria.fr/}}.
Furthermore, an important issue in the EMO field relates to the algorithm performance analysis and to set quality metrics~\cite{ZT+:03}.
A couple of metrics are featured within ParadisEO-MOEO, 
including the hypervolume metric in both its unary~\cite{ZT:99} and its binary~\cite{ZT+:03} form, 
the entropy metric~\cite{BST:02}, the contribution metric~\cite{MTR:00} as well as 
the additive and the multiplicative $\epsilon$-indicators~\cite{ZT+:03}.
Another interesting feature is the possibility to compare the current archive with the archive of the previous generation by means of a binary metric, 
and to print the progression of this measure iteration after iteration.

\subsubsection{EMO Algorithms}
\begin{figure}[htbp]
	\centering
		\includegraphics[scale=0.24]{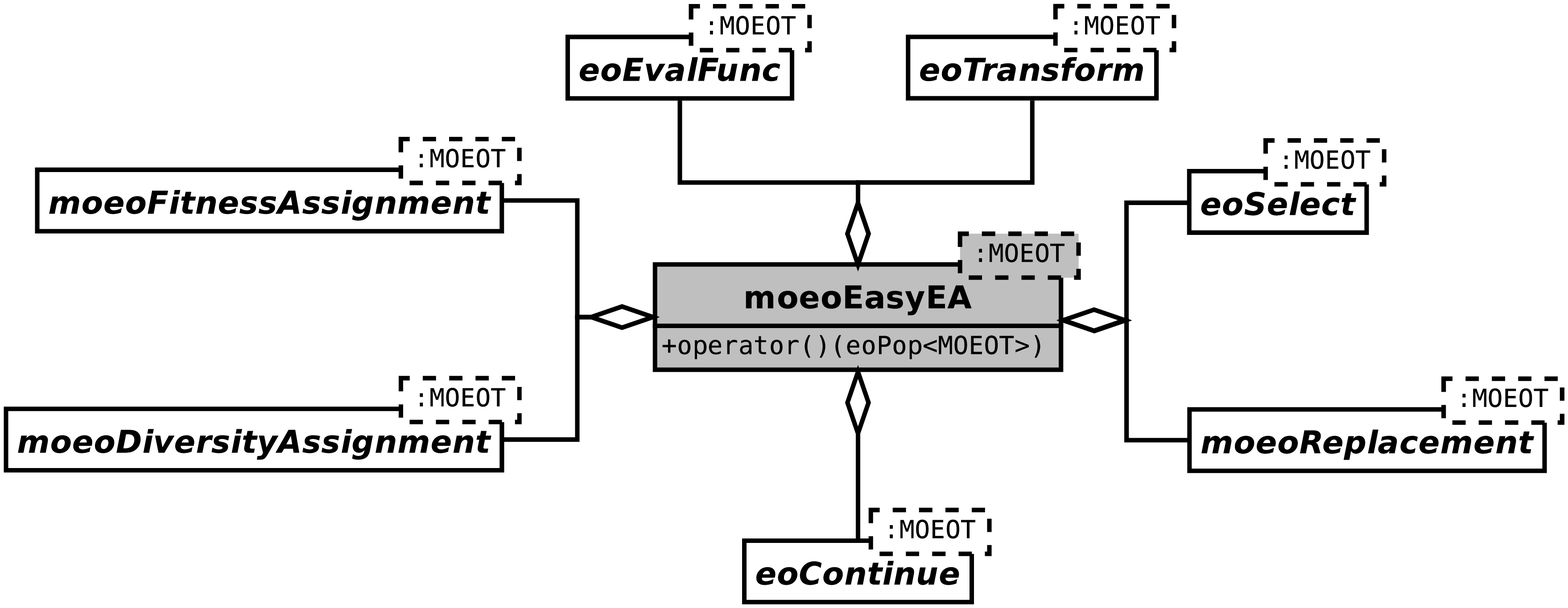}
	\caption{A possible instantiation for the design of an EMO algorithm.}
	\label{fig:uml:sketch}
\end{figure}

Now that all the basic, EA-related and EMO-specific components are defined, an EMO algorithm can easily be designed using the fine-grained classes of ParadisEO.
As the implementation is conceptually divided into components, 
different operators can be experimented without engendering significant modifications in terms of code writing.
As seen before, a wide range of components are already provided.
But, keep in mind that this list is not exhaustive as the framework perpetually evolves and offers all that is necessary to develop new ones with a minimum effort.
Indeed, ParadisEO is a white-box framework that tends to be flexible while being as user-friendly as possible. 
Fig.~\ref{fig:uml:sketch} illustrates the use of the \verb|moeoEasyEA| class that allows to define an EMO algorithm in a common fashion, 
by specifying all the particular components required for its implementation.
All classes use a template parameter {\itshape MOEOT (Multiobjective Evolving Object Type)} 
that defines the representation of a solution for the problem under consideration.
This representation might be implemented by inheriting of the \verb|MOEO| class as described in Sect.~\ref{sec:representation}.
Note that the archive-related components does not appear in the UML diagram, as we chose to let the use of an archive as optional.
The archive update can easily be integrated into the EA by means of the checkpointing process.
Similarly, the initialization process does not appear either, as an instance of \verb|moeoEasyEA| starts with an already initialized population.

\begin{figure}[htbp]
	\centering
		\includegraphics[scale=0.24]{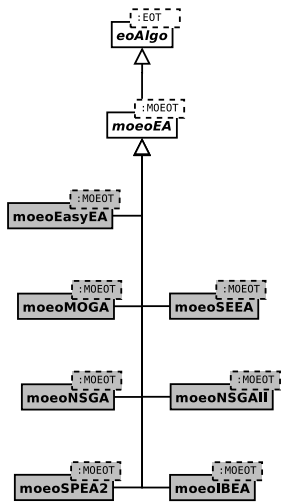}
	\caption{UML diagram for easy-to-use EMO algorithms.}
	\label{fig:uml:algo}
\end{figure}

In order to satisfy both the common user and the more experimented one, ParadisEO-MOEO also provides even more easy-to-use EMO algorithms, see Fig.~\ref{fig:uml:algo}.
These classes propose different implementations of some state-of-the-art methods by using the fine-grained components of ParadisEO.
They are based on a simple combination of components, as described in Sect.~\ref{sec:instances}.
Hence, MOGA~\cite{FF:93}, NSGA~\cite{SD:94}, NSGA-II~\cite{DA+:02}, SPEA2~\cite{ZLT:01}, IBEA~\cite{ZK:04} and SEEA~\cite{LJT:08}
are proposed in a way that a minimum number of problem-~or algorithm-specific parameters are required. 
For instance, to instantiate NSGA-II for a new continuous MOP, it is possible to use standard operators for representation, initialization and variation, 
so that the evaluation is the single component to be implemented.
These easy-to-use algorithms also tends to be used as references for a fair performance comparison in the academic world, 
even if they are also well-suited for a straight use to solve real-world MOPs.
In a near future, other easy-to-use EMO metaheuristics will be proposed while new fined-grained components will be implemented into ParadisEO-MOEO.

\subsection{Discussion}
ParadisEO-MOEO has been used and experimented to solve a large range of MOPs from both academic and real-world fields, what validates its high flexibility.
Indeed, various academic MOPs have been tackled within ParadisEO-MOEO, 
including continuous test functions (like the ZDT and DTLZ functions family defined in~\cite{DT+:05}), 
scheduling problems (permutation flow-shop scheduling problem~\cite{LB+:07b}), 
routing problems (multiobjective traveling salesman problem, bi-objective ring star problem~\cite{LJT:08}), and so on.
Moreover, it has been successfully employed to solve real-world applications in structural biology~\cite{BJ+:08}, 
feature selection in cancer classification~\cite{TJ+:08}, materials design in chemistry~\cite{SJ+:08}, etc. 
Besides, a detailed documentation as well as some tutorial lessons and problem-specific implementations are freely available on the ParadisEO website%
\footnote{\url{http://paradiseo.gforge.inria.fr}}.
And we expect the number of MOP contributions to largely grow in a near future.
Furthermore, note that the implementation of EMO algorithms is just an aspect of the features provided by ParadisEO.
Hence, hybrid mechanisms can be exploited in a natural way to make cooperating metaheuristics belonging to the same or to different classes.
Moreover, the three main parallel models are concerned: algorithmic-level, iteration-level and solution-level and are portable on different types of architecture.
Indeed, the whole framework allows to conveniently design hybrid as well as parallel and distributed metaheuristics, including EMO methods.
For instance, in the frame of ParadisEO, hybrid EMO algorithms have been experimented in~\cite{LJT:08},
a multiobjective cooperative island model has been designed in~\cite{TCM:07}, and costly evaluation functions have been parallelized in~\cite{BJ+:08}.
The reader is referred to~\cite{CMT:04} for more information about ParadisEO hybrid and parallel models.

\section{Concluding Remarks}
The current paper presents two complementary contributions: 
the formulation of a unified view for evolutionary multiobjective optimization and the description of a software framework for the development of such algorithms.
First, we identified the common concepts shared by many evolutionary multiobjective optimization techniques,
separating the problem-specific part from the invariant part involved in this class of resolution methods.
We emphasized the main issues of fitness assignment, diversity preservation and elitism.
Therefore, we proposed a unified conceptual model, based on a fine-grained decomposition,
and we illustrated its robustness and its reliability by treating a number of state-of-the-art algorithms as simple instances of the model.
Next, this unified view has been used as a starting point for the design and the implementation of a general-purpose software package called ParadisEO-MOEO.
ParadisEO-MOEO is a free C++ white-box object-oriented framework dedicated to the flexible and reusable design of evolutionary multiobjective optimization algorithms.
It is based on a clear conceptual separation between the resolution methods and the problem they are intended to solve,
thus conferring a maximum code and design reuse.
This global framework has been experimentally validated by solving a comprehensive number of both academic and real-world multiobjective optimization problems.

However, we believe that a large number of components involved in evolutionary multiobjective optimization are shared by many other search techniques.
Thereafter, we plan to generalize the unified model proposed in this paper to other existing metaheuristic approaches for multiobjective optimization.
Hence, multiobjective local search or scatter search methods might be interesting extensions to explore 
in order to investigate their ability and their modularity for providing such a flexible model as the one presented in this paper.
Afterwards, the resulting general-purpose models and their particular mechanisms would be integrated into the ParadisEO-MOEO software framework.

\section*{Acknowledgment}
The authors would like to gratefully acknowledge Thomas Legrand, J\'er\'emie Humeau, and Abdel-Hakim Deneche 
for their helpful contribution on the implementation part of this work,
as well as S\'ebastien Cahon and Nouredine Melab for their work on the preliminary version of the ParadisEO-MOEO software framework presented in this paper.
This work was supported by the ANR DOCK project.

\bibliographystyle{IEEEtran}
\bibliography{IEEEabrv,mo,meta,frameworks,appli}

\newpage
\tableofcontents

\end{document}